\documentclass[12pt]{amsart}
\usepackage{amssymb}
\usepackage{amsmath}
\usepackage{mathabx}
\usepackage{amsthm}
\usepackage{amstext}
\usepackage{amsopn}
\usepackage{mathrsfs}
\usepackage{latexsym}
\DeclareMathAlphabet{\mathpzc}{OT1}{pzc}{m}{it}
\allowdisplaybreaks
\newtheorem{Definition}{Definition}[section]

\newtheorem{Theorem}{Theorem}[section]
\newtheorem{Corollary}{Corollary}[section]
\newtheorem{Remark}{Remark}[section]
\newtheorem{Example}{Example}[section]

\begin{document}
\bibliographystyle{plain}
\title[Limit distributions of Gaussian block ensembles]
{Limit distributions of Gaussian block ensembles}
\author[R. Lenczewski]{Romuald Lenczewski}
\address{Romuald Lenczewski, \newline
Instytut Matematyki i Informatyki, Politechnika Wroc\l{}awska, \newline
Wybrze\.{z}e Wyspia\'{n}skiego 27, 50-370 Wroc{\l}aw, Poland  \vspace{10pt}}
\email{Romuald.Lenczewski@pwr.wroc.pl}
\begin{abstract}
It has been shown by Voiculescu that important classes of square independent random matrices 
are asymptotically free, where freeness is a noncommutative analog of classical independence.
Recently, we introduced the concept of matricial freeness, which is similar to freeness in free probability, 
but it also has some matricial features. Using this new concept of noncommutative independence, we 
described the asymptotics of blocks and symmetric blocks of certain classes of independent random matrices.
In this paper, we present the main results obtained in this framework, concentrating on the ensembles of blocks of Gaussian random matrices.
\end{abstract}

\maketitle

\section{Introduction}
Our main objective is to describe asymptotic joint distributions of rectangular blocks of 
independent random matrices, called {\it random blocks}, under the expectation of normalized partial traces. 
For that purpose, we use a new concept of noncommutative independence 
called {\it matricial freeness} and associated arrays of operators which give Hilbert space 
realizations of these distributions. 
The well-known connection between free probability and the asymptotics 
of independent square random matrices under the expectation of normalized trace 
can also be reproduced in this framework.

The most fundamental results of this nature were obtained by Voiculescu [16], 
who showed that certain ensembles of independent $n\times n$ random matrices $\{Y(u,n):u\in \mathpzc{U}\}$ 
were asymptotically free under the expectation of the normalized complete trace, 
$\tau(n)=\mathbb{E}\otimes {\rm Tr}(n)$, where ${\rm Tr}(n)$ stands for the trace divided by $n$.
In particular, if the entries of $Y(u,n)$ are i.i.d. complex $(0,1/n)$-Gaussian random variables, 
we can symbolically write
$$
\lim_{n\rightarrow \infty}Y(u,n)=\eta(u)
$$
where $\{\eta(u):u\in \mathpzc{U}\}$ is the standard free circular system of operators and convergence is understood 
in the sense of mixed moments under $\tau(n)$. The operators $\eta(u)$ live in the free Fock space 
and have the standard circular distribution (uniform distribution on the unit disc in the complex plane) 
in the vacuum state. A similar result holds for Hermitian Gaussian random matrices 
whose limit joint distributions are described by mixed moments of 
free Gaussian operators with semicircle distributions.

Gaussian random matrices studied by Voiculescu had i.i.d. entries, except that in the case of 
Hermitian ensembles it holds that $\,\overline{Y_{i,j}(u,n)}=Y_{j,i}(u,n)$. 
If the entries are independent but not identically distributed, 
standard free probability may not suffice to describe asymptotic distributions of these matrices.  
One of the approaches is to employ the much more general scheme of freeness with 
amalgamation as in the papers of Shlyakhtenko [14] and Benaych-Georges [3]. 
Recently, we have developed another approach [7,8] in the case when the entries are independent and block-identically distributed
({\it i.b.i.d.}). This approach is based on the concept of matricial freeness [6], 
which can be viewed as a matricial generalization of freeness, lying somewhere between freeness and freeness with amalgamation.

In this context, we decompose matrices $Y(u,n)$ with i.b.i.d. entries 
into rectangular blocks $\{S_{p,q}(u,n):1\leq p,q \leq r\}$, or symmetric blocks $\{T_{p,q}(u,n):1\leq p \leq q \leq r\}$. The block analog of the above Voiculescu's result for the ensemble of 
independent Gaussian matrices can then be written in the form
$$
\lim_{n\rightarrow \infty}S_{p,q}(u,n)=\zeta_{p,q}(u)
$$
where $\zeta_{p,q}(u)$ are certain operators living in the 
{\it matricially free Fock space}, a matricial analog of the free Fock space, into which $\eta(u)$ decompose.
The convergence is in the sense of mixed *-moments under the expectation of normalized partial traces 
over the subsets of basis vectors defined by the block decomposition, denoted 
$\tau_{q}(n)={\mathbb E}\circ {\rm Tr}_{q}(n)$, where $1\leq q\leq r$.
In turn, the corresponding moments of operators are computed with respect to the family of vacuum states. 
A similar theorem holds for symmetric blocks of Hermitian and non-Hermitian Gaussian random matrices.

This block model encodes more degrees of freedom than the usual framework of square 
random matrices with i.i.d. Gaussian entries. Apart from block variances, the most important parameters 
are the {\it asymptotic dimensions} of blocks, namely
$$
d_{q}=\lim_{n\rightarrow \infty}\frac{n_{q}}{n}, 
$$
where $n=n_1+\ldots +n_{r}$ is the $n$-dependent partition of $n$ defined by the block decomposition of 
$Y(u,n)$. We assume that these limits exist, but it is possible that some of them vanish. This leads 
to three types of blocks or symmetric blocks: {\it balanced} (if both their asymptotic dimensions are positive), 
{\it unbalanced} (if one asymptotic dimension is positive and the other one vanishes) 
and {\it evanescent} (if both asymptotic dimensions vanish). For details, see [8].

Using all types of blocks, we can construct random matrix models for other notions of independence, 
such as monotone independence, boolean independence, c-freeness (conditional freeness) and s-freeness (freeness with subordination). 
Moreover, the limit moments can be written in a quite explicit form as polynomials in $d_1, \ldots , d_r$, with 
block variances as additional parameters. In some cases, this enables us to construct simple random matrix models for certain families of probability measures, like the model for free Meixner laws [1,4] constructed in [9].
At the same time, studying such polynomials can lead to new results in free probability. 
For instance, a study of the limit distributions of products of independent rectangular Gaussian random matrices 
produced polynomials which can be viewed as {\it multivariate Fuss-Narayana polynomials} and, moreover, turned out to be the moments of the free multiplicative convolution of Marchenko-Pastur distributions with arbitrary shape parameters [10].

In order to make this paper easy to follow, we avoid technical details and only sketch some proofs. 
For complete proofs, we refer the reader to [7,8,9,10].

\section{Matricial operator systems}

In the framework of free probability, an important role is played by semicircular and circular systems of operators. Such systems were introduced by Voiculescu, who used them to construct the generators of free group factors in order to prove certain isomorphisms between them [16].

These generators play also an important role in the framework of matricial freeness. 
For the definition of that concept, which can be viewed as a matricial generalization of freeness
involving a family of states instead of one state, we refer the reader to [6]. 
In this paper, we prefer to restrict our attention to important examples of arrays of operators [7,8] which are matricially free with respect to a family of states associated with vacuum vectors in the `matricial' analog of the free Fock space, called the matricially free Fock space.

\begin{Definition}
{\rm  Let $({\mathcal H}_{p,q}(u))$, $u \in \mathpzc{U}$, be a family
of $r\times r$ arrays of Hilbert spaces. The {\it matricially free Fock space} is the 
Hilbert space direct sum of the form
\begin{equation*}
{\mathcal M}= \bigoplus_{q=1}^{r} {\mathcal M}_{q},
\end{equation*}
where 
\begin{equation*}
{\mathcal M}_{q}={\mathbb C}\Omega_{q}\oplus \bigoplus_{m=1}^{\infty}
\bigoplus_{\stackrel{p_1,\ldots , p_m\in [r]}
{\scriptscriptstyle u_1, \ldots , u_n\in \mathpzc{U}}
}
{\mathcal H}_{p_1,p_2}(u_{1})\otimes {\mathcal H}_{p_2,p_3}(u_2)
\otimes \ldots \otimes 
{\mathcal H}_{p_m,q}(u_{m}),
\end{equation*}
with $[r]:=\{1, 2, \ldots, r\}$ and $\{\Omega_1, \ldots , \Omega_{r}\}$ being the 
set of unit (vacuum) vectors, equipped with the canonical inner product. 
Let $\{\Psi_1, \ldots , \Psi_r\}$ be the corresponding 
states replacing the single vacuum state in the free Fock space.}
\end{Definition}

In fact, for most purposes, it suffices to take each Hilbert space to be one-dimensional, namely
${\mathcal H}_{p,q}(u)={\mathbb C}e_{p,q}(u)$, where the vectors $e_{p,q}(u)$ 
form an orthonormal basis. Arrays of certain partial isometries living in this Fock 
space which remind free creation operators $\ell(u)$ living in the free Fock space serve as generators of 
certain Toeplitz-Cuntz-Krieger algebras [5]. These partial isometries, 
when multiplied by positive scalars, become matricial analogs of free creation operators.

\begin{Definition}
{\rm Let $B(u)=(b_{p,q}(u))$ be an array of positive real numbers for any $u\in \mathpzc{U}$.
We associate with each such array the array of {\it matricially free creation operators} whose 
non-trivial action onto the basis vectors is  
\begin{eqnarray*}
\wp_{p,q}(u)\Omega_q&=&\sqrt{b_{p,q}(u)}e_{p,q}(u)\\
\wp_{p,q}(u)(e_{q,t}(s))&=&\sqrt{b_{p,q}(u)}(e_{p,q}(u)\otimes e_{q,t}(s))\\
\wp_{p,q}(u)(e_{q,t}(s)\otimes w)&=&\sqrt{b_{p,q}(u)}(e_{p,q}(u)\otimes e_{q,t}(s)\otimes w)
\end{eqnarray*}
for any $p,q,t\in [r]$ and $u,s\in \mathpzc{U}$, where $e_{q,t}(s) \otimes w$ is a basis vector. Their actions onto the remaining basis vectors give zero. 
The corresponding {\it matricially free annihilation operators} are their adjoints 
$\wp_{p,q}^{*}(u)$. In some cases, it will be convenient to use the same notation even if $b_{p,q}(u)=0$, when 
we obtain trivial operators.}
\end{Definition}

Certain linear combinations of the matricially free creation and annihilation operators are of special interest.
We have studied several matricial systems of operators constructed from the matricially free creation and annihilation operators 
or their symmetrized counterparts when describing the asymptotic distributions of blocks or symmetric blocks of Hermitian and non-Hermitian Gaussian random matrices [7,8]. We define them below, mentioning also the corresponding ensembles of random blocks. We also give their realizations as operator-valued matrices, which identifies 
some of them with the generators of free group factors used by Voiculescu in [16].
\begin{enumerate}
\item
Non-trivial matricially free creation and annihilation operators can be 
realized as operator-valued matrices from the $C^*$-algebra ${\mathcal A}\otimes M_{n}({\mathbb C})$, 
where ${\mathcal A}$ is the $C^*$-algebra generated by the family 
$\{\ell(p,q,u):1\leq p,q \leq r, u\in \mathpzc{U}\}$ of *-free creation operators living in the free Fock space. Namely, 
\begin{eqnarray*}
\wp_{p,q}(u)&=&\ell(p,q,u)\otimes e(p,q)\\
\wp_{p,q}^*(u)&=&\ell(p,q,u)^*\otimes e(q,p),
\end{eqnarray*}
where $\{e(p,q): 1\leq p,q \leq r\}$ is the array of matrix units.
If $\varphi$ is the vacuum state in this free Fock space and 
$\psi_q$ is the vector state on $M_{n}({\mathbb C})$ defined by the basis vector $e_{q}$ of ${\mathbb C}^{n}$, then the state $\Psi_q$
can be identified with $\varphi\otimes \psi_q$, as we showed in 
[11].
\item
{\it Matricially free Gaussian operators} are self-adjoint operators of the form
$$
\omega_{p,q}(u)=\wp_{p,q}(u)+\wp_{p,q}(u)^*,
$$
and they are the canonical Gaussian operators in our framework. It turns out that 
they describe the limit distributions of unbalanced symmetric blocks of Hermitian Gaussian random matrices.
\item
It is convenient to introduce the {\it symmetrized creation operators} as
$$
\widehat{\wp}_{p,q}(u)=
\left\{\begin{array}{ll}
\wp_{p,q}(u)+\wp_{q,p}(u) & {\rm if} \; p<q\\
\wp_{q,q}(u)& {\rm if}\; p=q
\end{array}
\right.
$$
and the {\it symmetrized annihilation operators} as their adjoints $\widehat{\wp}^*_{p,q}(u)$. 
\item
In order to describe the limit distributions of balanced symmetric blocks of Hermitian Gaussian random matrices, we need to use 
{\it symmetrized Gaussian operators} 
$$
\widehat{\omega}_{p,q}(u)=\widehat{\wp}_{p,q}(u)+\widehat{\wp}_{p,q}(u)^{*},
$$
where all operators are non-trivial.
It is easy to see that if all creation and annihilation operators involved are non-trivial, then 
the above ones can be identified with the operator-valued matrices of the form
$$
\widehat{\omega}_{p,q}(u)=\left\{
\begin{array}{ll}
g(p,q,u)\otimes e(p,q)+g(p,q,u)^*\otimes e(q,p) & {\rm if}\; p<q\\
s(q,u)\otimes e(q,q)& {\rm if} \;p=q
\end{array}
\right.
$$
where $\{g(p,q,u): 1\leq p \leq q \leq r, u\in \mathpzc{U}\}$ is a family of 
circular operators and $\{s(q,u): 1\leq q \leq r, u\in \mathpzc{U}\}$ 
is a family of semicircular operators. These matrices are generators of free group factors
introduced by Voiculescu [16].
\item
In order to describe the limit distributions of the usual (non-symmetric) blocks of non-Hermitian 
Gaussian random matrices, we need to use {\it matricial $R$-circular operators} [11]
$$
\zeta_{p,q}(u)=
\wp_{p,q}(u')+\wp_{q,p}(u'')^*,
$$
where the notation means that in order to construct one such operator labelled by $u$ one needs two 
operators, which are labelled by $u'$ and $u''$, thus the index set $\mathpzc{U}$ has to be
doubled. If the added operators are non-trivial, then we can identify the above ones
with operator-valued matrices of the form
$$
\zeta_{p,q}(u)=g(p,q,u)\otimes e(p,q)
$$
for any $1\leq p,q\leq r$ and any $u\in \mathpzc{U}$, where 
$\{g(p,q,u): 1\leq p, q \leq r, u\in \mathpzc{U}\}$ is a family of circular operators (here,
the covariances of $\wp_{p,q}(u)$ are symmetric and identical for $u'$ and $u''$). Sums $\zeta(u)=\sum_{p,q}\zeta_{p,q}(u)$
are $R$-cyclic matrices introduced in [13].
\item
In order to describe the limit joint distributions of the symmetric blocks of 
non-Hermitian Gaussian random matrices, we need to use {\it matricial circular operators} [8]
which are natural symmetrizations of $\zeta_{p,q}(u)$, namely
$$
\eta_{p,q}(u)=\widehat{\wp}_{p,q}(u')+\widehat{\wp}_{p,q}(u'')^*,
$$
where $u\in \mathpzc{U}$ and $1\leq p\leq q\leq r$. If the added operators are non-trivial, then
the above ones can be identified with the operator-valued matrices of the form
$$
\eta_{p,q}(u)=\left\{
\begin{array}{ll}
g(p,q,u)\otimes e(p,q)+g(q,p,u)\otimes e(q,p) & {\rm if}\; p<q\\
s(q,u)\otimes e(q,q)& {\rm if} \;p=q
\end{array}
\right.
$$
where $\{g(p,q,u): 1\leq p,q \leq r, u\in \mathpzc{U}\}$ is a family of 
circular operators and $\{s(q,u): 1\leq q \leq r, u\in \mathpzc{U}\}$ 
is a family of semicircular operators.  
\end{enumerate}

\section{Hermitian Gaussian Block Ensemble}

In free probability, we study the mixed moments of independent $n\times n$  random matrices 
$\{Y(u,n): u\in \mathpzc{U}\}$ under the expectation of the normalized trace 
$$
\tau(n)={\mathbb E}\circ {\rm Tr}(n)
$$ 
where 
$$
{\rm Tr}(n)(A)=\frac{1}{n}\, {\rm Tr}(A)
$$ 
for any $n\times n$ matrix $A$ as $n\rightarrow \infty$. Let us recall the asymptotic freeness result of Voiculescu for independent Hermitian Gaussian random matrices [15].

\begin{Theorem}
If we are given an ensemble of independent Hermitian $n\times n$ random matrices 
$\{Y(u,n):u\in \mathpzc{U}\}$, whose entries $Y_{i,j}(u,n)$ 
satisfy $Y_{i,j}(u,n)=\overline{Y_{j,i}(u,n)}$, 
are complex $(0,1/n)$-Gaussian if $i\neq j$ and real $(0,1/n)$-Gaussian if $i=j$ and 
the family $\{Y_{i,j}(u,n):1\leq i \leq j \leq n\}$ is independent for any $u$, then 
$$
\lim_{n\rightarrow \infty}Y(u,n)=
\omega(u),
$$
which should be understood as the convergence of mixed moments of matrices 
under $\tau(n)$ to the mixed moments of the corresponding free standard 
semicircular operators $\omega(u)$ under the vacuum state $\Phi$ in the free Fock space.
\end{Theorem}

\begin{Remark}
{\rm Let us make some remarks which will enable us to formulate our main results. 
\begin{enumerate}
\item
We stated the above theorem in a simplified form, which will also be used 
when we state our results on the asymptotic distributions of random blocks. A more explicit 
formulation says that 
$$
\lim_{n\rightarrow \infty}\tau(n)(Y(u_1,n)\ldots Y(u_m,n))=
\Phi(\omega(u_1)\ldots \omega(u_m)),
$$ 
for any $u_1,\ldots , u_m\in \mathpzc{U}$, where $\{\omega(u): u\in \mathpzc{U}\}$
is a standard free semicircular family, which means in particular 
that $\Phi(\omega(u)^{2})=1$ for any $u$. Alternatively, one could say that
the family $\{Y(u,n):u\in\mathpzc{U}\}$ is asymptotically free under $\tau(n)$ 
and that the limit distribution of each $Y(u,n)$ is the standard semicircular 
Wigner distribution $W(0,1)$ with density $(2\pi)^{-1}\sqrt{4-x^2}$ on $[-2,2]$.
\item
In order to define a family of partial traces, take the decomposition of the set $[n]:=\{1,2, \ldots , n\}$
into $r$ disjoint intervals
$$
[n]:=N_1\cup \ldots \cup N_r
$$ 
and denote by $n_q$ the cardinality of $N_{q}$. Next, let 
$$
I(n)=D_{1}+\ldots +D_q
$$ 
be the corresponding decomposition of the $n\times n$  identity matrix $I(n)$, that is
$(D_k)_{i,j}=1$ whenever $i=j\in N_k$, with the remaining entries equal to zero, 
where $1\leq k \leq r$. Of course, the objects $N_q,n_q,D_q$ depend on $n$, but this is supressed in the 
notation.
\item
By {\it partial traces} we then understand states of the form
$$
\tau_{q}(n)={\mathbb E}\circ {\rm Tr}_{q}(n)
$$ 
where 
$$
{\rm Tr}_q(n)(A)=\frac{1}{n_q}\, 
{\rm Tr}(D_qAD_q)
$$ 
and $1\leq q \leq r$. 
\item
By {\it random blocks} of a random matrix $Y(u,n)$ we shall understand $n\times n$ matrices of the form
$$
S_{p,q}(u,n)=D_{p}Y(u,n)D_{q},
$$
for any $1\leq p,q\leq r$ and $u\in \mathpzc{U}$, where $r,n\in {\mathbb N}$. 
In particular, if $Y(u,n)$ is Hermitian, then $S_{q,p}(u,n)=S_{p,q}(u,n)^*$.
Clearly, we have the decomposition
$$
Y(u,n)=\sum_{1\leq p,q\leq r}S_{p,q}(u,n)
$$
for any $u,n$.
\item
We assume that all variables $Y_{i,j}(u,n)$ which belong to the same block $S_{p,q}(u,n)$ have the same 
covariance $\mathbb E(\overline{Y_{i,j}(u,n)}Y_{i,j}(u,n))=v_{p,q}(u)/n$ whenever
$(i,j)\in N_p\times N_q$. We denote by $V(u)=(v_{p,q}(u))$ the corresponding matrices of block covariances. 
Apart from the dimension matrix 
$$
D={\rm diag}(d_1, \ldots , d_r),
$$
these matrices are additional parameters of the ensemble.  In fact, the limit distributions will 
be expressed in terms of matrices $B(u)=DV(u)$.
\item
By {\it random symmetric blocks} of a random matrix $Y(u,n)$ we shall understand $n\times n$ matrices of the form 
$$
T_{p,q}(u,n)=
\left\{
\begin{array}{ll}
S_{p,q}(u,n)+S_{q,p}(u,n) & {\rm if}\;p< q\\
S_{q,q}(u,n) & {\rm if} \;p=q
\end{array}
\right.
$$
for any $1\leq p \leq q\leq r$ and $u\in \mathpzc{U}$, where $r,n\in {\mathbb N}$.
We have the decomposition 
$$
Y(u,n)=\sum_{1\leq p \leq q \leq r}T_{p,q}(u,n)
$$
for any $u,n$. Clearly, if $Y(u,n)$ is Hermitian, then 
$T_{p,q}(u,n)=T_{p,q}(u,n)^*$.
\item
When speaking of limit distributions of mixed moments of symmetric blocks, we will 
use a simplified formulation, similar to that in the free case. For instance 
$$
\lim_{n\rightarrow \infty}T_{p,q}(u,n)=\widehat{\omega}_{p,q}(u)
$$
in the Hermitian case will mean that
$$
\lim_{n\rightarrow \infty}
\tau_q(n)(T_{p_1,q_1}(u_1,n)\ldots T_{p_m,q_m}(u_m,n))
$$
$$
=
\Psi_q(\widehat{\omega}_{p_1,q_1}(u_1)\ldots \widehat{\omega}_{p_m,q_m}(u_m))
$$
for any $1\leq p_1\leq q_1\leq r, \ldots , 1\leq p_m\leq q_m\leq r$, $1\leq q\leq r$ 
and $u_1, \ldots , u_m\in\mathpzc{U}$. The operators will always belong to one of the 
families defined in Section 2. A similar formulation will be used for blocks $S_{p,q}(u,n)$.
\end{enumerate}
}
\end{Remark}

We can describe limit joint distributions of blocks and symmetric blocks of Hermitian random matrices 
under rather general assumptions. It suffices to assume that the family $\{Y(u,n):u\in \mathpzc{U}\}$ is asymptotically 
free and asymptotically free from $\{D_1, \ldots , D_r\}$ under $\tau(n)$ and that their norms are
uniformly boubded almost surely. This class includes unitarily invariant random matrices whose limit
moments are compactly supported probability measures on the real line. This general version was proved in [8, Theorem 6.1].
For the sake of simplicity, we restrict our attention here to the case of Hermitian Gaussian random matrices.

\begin{Theorem}
If we are given an ensemble $\{Y(u,n):u\in \mathpzc{U}\}$ of independent Hermitian $n\times n$ random matrices
whose entries $Y_{i,j}(u,n)$ satisfy $Y_{i,j}(u,n)=\overline{Y_{j,i}(u,n)}$, are complex $(0,v_{p,q}(u)/n)$-Gaussian if $i\neq j$ and 
$(i,j)\in N_p\times N_q$, and real $(0,v_{q,q}(u)/n)$-Gaussian if $i=j\in N_q$ and the family $\{Y_{i,j}(u,n):1\leq i \leq j \leq n\}$ 
is independent for any $u$, then 
$$
\lim_{n\rightarrow \infty}T_{p,q}(u,n)=
\widehat{\omega}_{p,q}(u)
$$
for any $p\leq q$ and $u\in \mathpzc{U}$,
where convergence is in the sense of mixed moments under partial traces and the arrays $(\widehat{\omega}_{p,q}(u))$ 
are associated with the symmetric covariance matrices $B(u)=DV(u)$, respectively.
\end{Theorem}
{\it Sketch of the proof.}
First, let us assume that $v_{p,q}(u)=1$ for any $p,q,u$. In that case, we can use asymptotic freeness of the 
family $\{Y(u,n):u\in \mathpzc{U}\}$ and its asymptotic freeness with respect to the family $\{D_1, \ldots, D_r\}$
(these are deterministic diagonal matrices) to describe the limit joint distributions of the blocks $S_{p,q}(u,n)$
under $\tau(n)$ since
$$
S_{p,q}(u,n)=D_{p}Y(u,n)D_{q}
$$
for any $p,q,u$. Moreover, we know that the limit distribution of each 
$Y(u,n)$ is the standard semicircle Wigner law $W(0,1)$ and a direct computation 
gives $\tau(n)(D_{q})=n_q/n\rightarrow d_{q}$
for any $q$. Now, it is not hard to show that the family of operators
$$
\omega(u):=\sum_{p,q}\omega_{p,q}(u),
$$
where $\omega_{p,q}(u)=\wp_{p,q}(u)+\wp_{p,q}(u)^*$ 
for any $u$ and $\wp_{p,q}(u)$ has covariance $d_p$, 
is free with respect to $\Psi=\sum_{q}d_{q}\Psi_q$. 
We abuse the notation a little since by $\omega(u)$ we denoted 
a standard semicircular operator on the free Fock space, but this is justified by the fact
that each $\omega(u)$ has the standard semicircle distribution under $\Psi$. 
Therefore, we have in fact a realization of the standard free semicircular family 
$\{\omega(u): u\in \mathpzc{U}\}$ as operators living in the matricially free Fock space ${\mathcal M}$. 
Now, it suffices to find an appropiate limit realization 
for $\{D_1, \ldots , D_r\}$. We have shown in [8] that it is given by the family $\{P_1, \ldots , P_r\}$, where
$$
P_{q}=1\otimes e(q,q),
$$
using the tensor product realization described in Section 2, since the family
$\{\omega(u): u\in \mathpzc{U}\}$
is free from $\{P_1, \ldots , P_r\}$ and $\Psi(P_q)=d_q$. Therefore, any limit mixed *-moment 
of the random blocks $S_{p,q}(u)$ under $\tau(n)$ can be expressed as a mixed
*-moment of the corresponding operators $P_{p}\,\omega(u)P_q$ under $\Psi$. Symbolically, 
since $\lim_{n\rightarrow \infty}Y(u,n)=\omega(u)$, we have
$$
\lim_{n\rightarrow \infty}S_{p,q}(u,n)=P_{p}\,\omega(u)P_{q}=\wp_{p,q}(u)+\wp_{q,p}(u)^*:=\varsigma_{p,q}(u)
$$
except that the moments of blocks are computed under $\tau(n)$ and those of the operators $P_{p}\,\omega(u)P_q$ 
are computed under $\Psi$. In order to pass from $\tau(n)$ to $\tau_q(n)$, notice that
$$
\tau_q(n)(A)=\frac{n}{n_{q}}\tau (n)(D_{q}AD_q)
$$
so the limit mixed *-moments of random blocks under $\tau_q(n)$ 
are given by the mixed *-moments of the above operators under $\Psi_q$, where $1\leq q \leq r$. All this holds
provided $d_{q}\neq 0$. The case when $d_q=0$ is slightly more complicated and is omitted 
here (we refer the reader to [8]). It easily follows that 
$$
\lim_{n\rightarrow \infty}T_{p,q}(u,n)=\widehat{\omega}_{p,q}(u)
$$
for any $1\leq p\leq q \leq r$ and $u\in \mathpzc{U}$.
It can also be seen that one can rescale blocks $S_{p,q}(u,n)$, which means that one can rescale
the block covariances and take arbitrary non-negative $v_{p,q}(u)$ (except that we must have the symmetry $v_{p,q}(u)=v_{q,p}(u)$ 
since $Y(u,n)$ is Hermitian). This proves that in the case when the covariances are equal to $v_{p,q}(u)$
within block $S_{p,q}(u)$, the limit operator $\varsigma_{p,q}(u)$ gets rescaled by $v_{p,q}(u)$ and thus 
the covariance of $\wp_{p,q}(u)$ becomes $b_{p,q}(u)=d_{p}v_{p,q}(u)$. This completes the proof.
\hfill $\blacksquare$

\begin{Corollary}
In particular,
\begin{enumerate}
\item
if $d_p=1$ and $d_q=0$, then $\lim_{n\rightarrow \infty}T_{p,q}(u,n)=\omega_{p,q}(u)$,
\item
if $d_p=0$ and $d_q=1$, then $\lim_{n\rightarrow \infty}T_{p,q}(u,n)=\omega_{q,p}(u)$,
\item
if $d_p=0$ and $d_q=0$, then $\lim_{n\rightarrow \infty}T_{p,q}(u,n)=0$.
\end{enumerate}
\end{Corollary}
{\it Proof.}
If $d_p=1$ and $d_q=0$, then $\omega_{q,p}(u)=0$ since $b_{q,p}(u)=d_{q}v_{q,p}(u)=0$ and thus 
$\widehat{\omega}_{p,q}(u)$ reduces to $\omega_{p,q}(u)$. In turn, if 
$T_{p,q}(u,n)$ is evanescent, then $\omega_{p,q}(u)=\omega_{q,p}(u)=0$ and thus
$\widehat{\omega}_{p,q}(u)=0$, which completes the proof.
\hfill $\blacksquare$.

\begin{Remark}
{\rm 
It should be remarked that the matricially free Gaussian operators are not operatorial realizations of 
the usual (non-symmetric) blocks $S_{p,q}(u)$ of Gaussian Hermitian matrices. In fact, it follows from the proof of Theorem 3.2 that 
$$
\lim_{n\rightarrow \infty}S_{p,q}(u)=\left\{
\begin{array}{ll}
g(p,q,u)\otimes e(p,q)& {\rm if}\; p\neq q\\
s(q,u)\otimes e(q,q)& {\rm if}\;p=q
\end{array}
\right.
$$
for any $1\leq p,q\leq r $ and $u\in \mathpzc{U}$, where convergence is understood 
in the usual sense (mixed moments of blocks under partial traces $\tau_{q}(n)$ converge to mixed moments
of the corresponding operators under $\Psi_q$). Here, $\{g(p,q,u): 1\leq p \leq q \leq r, u\in \mathpzc{U}\}$ is a family of circular operators and $\{s(q,u): 1\leq q \leq r, u\in \mathpzc{U}\}$ 
is a family of semicircular operators and we assume that $g(q,p,u)=g(p,q,u)^*$ for $p<q$.}
\end{Remark}

\begin{Example}
{\rm Fix $u\in \mathpzc{U}$ (omitted in our notations) and $p\neq q$. Let $g(p,q)=\ell_1+\ell_2^*$
and thus $g(q,p)=\ell_1^*+\ell_2$ for some *-free creation operators 
$\ell_1, \ell_2$. In view of Remark 3.2, we get
\begin{eqnarray*}
\lim_{n\rightarrow \infty}\tau_{q}(n)(S_{q,p}(n)S_{p,q}(n)S_{q,p}(n)S_{p,q}(n))&=&
\Psi_q(\varsigma_{q,p}\varsigma_{p,q}\varsigma_{q,p}\varsigma_{p,q})\\
&=&
\varphi((\ell_1^*+\ell_2)(\ell_1+\ell_2^*)(\ell_1^*+\ell_2)(\ell_1+\ell_2^*))\\
&=&
(d_{p}^2+d_pd_q)v_{p,q}^2
\end{eqnarray*}
where $\varphi(\ell_1^*\ell_1)=\Psi_q(\wp_{p,q}^*\wp_{p,q})=d_{p}v_{p,q}$ and  
$\varphi(\ell_2^*\ell_2)=\Psi_p(\wp_{q,p}^*\wp_{q,p})=d_pv_{q,p}$ 
and we use the symmetry $v_{p,q}=v_{q,p}$.}
\end{Example}

\section{Ginibre Block Ensemble}

Blocks of non-Hermitian Gaussian random matrices can be treated in a similar way. 
The associated ensmble of blocks cna be called the {\it Ginibre Block Ensemble}. The main difference is that the limit joint distributions of blocks are described by matricial $R$-circular operators 
$\zeta_{p,q}(u)$ and those of the symmetric blocks by matricial circular operators $\eta_{p,q}(u)$. 

Let us recall Voiculescu's theorem on the asymptotic freeness of the ensemble of independent Gaussian 
random matrices with i.i.d. entries (Ginibre Ensemble) [15].

\begin{Theorem}
If we are given an ensemble of independent $n\times n$ random matrices 
$\{Y(u,n):u\in \mathpzc{U}\}$, whose entries $Y_{i,j}(u,n)$ are complex $(0,1/n)$-Gaussian for any $i,j$
and the family $\{Y_{i,j}(u,n):1\leq i, j \leq n\}$ is independent, then 
$$
\lim_{n\rightarrow \infty}Y(u,n)=
\eta(u),
$$
where convergence is in the sense of mixed moments of matrices 
under $\tau(n)$ to the mixed moments of the corresponding free standard 
circular operators $\eta(u)$ under the vacuum state $\Phi$ in the free Fock space.
\end{Theorem}

Let us now formulate an analogous theorem for blocks and symmetric blocks of Gaussian random matrices 
with i.b.i.d. entries.

\begin{Theorem}
If we are given an ensemble of independent $n\times n$ random matrices 
$\{Y(u,n):u\in \mathpzc{U}\}$, whose entries $Y_{i,j}(u,n)$ are complex $(0,v_{p,q}(u)/n)$-Gaussian for any $(i,j)\in N_p\times N_q$
and the family $\{Y_{i,j}(u,n):1\leq i,j \leq n\}$ is independent, then 
$$
\lim_{n\rightarrow \infty}S_{p,q}(u,n)=
\zeta_{p,q}(u),
$$
where convergence is in the sense of mixed moments of blocks 
under partial traces $\tau_q(n)$ to the mixed moments of the corresponding 
matricial systems of operators under the vacuum states $\Psi_q$, respectively,
in the matricially free Fock space.
\end{Theorem}
{\it Sketch of the proof.}
The proof is similar to that of Theorem 3.1. In the case when all variables are i.i.d., 
we can use the asymptotic *-freeness of the family $\{Y(u,n):u\in \mathpzc{U}\}$ under $\tau(n)$ as $n\rightarrow \infty$ as well as 
their asymptotic *-freeness from the family of diagonal matrices. 
One can check that the family $\{\eta(u): u\in \mathpzc{U}\}$, where
$$
\eta(u)=\sum_{p,q}\zeta_{p,q}(u)=\sum_{p,q}g(p,q,u)\otimes e(p,q),
$$
is *-free under $\Psi$ as well as  *-free with respect to $\{P_1, \ldots , P_r\}$ under $\Psi$, 
where again $P_q=1\otimes e(q,q)$. Moreover, each $\eta(u)$ has the standard circular distribution under $\Psi=\sum_{q}\Psi_q$.
Since the asymptotic joint distribution of $\{D_1, \ldots , D_r\}$ under $\tau(n)$ agrees with that of 
$\{P_1, \ldots, P_r\}$ under $\Psi$, we must have
$$
\lim_{n\rightarrow \infty}S_{p,q}(u,n)=\lim_{n\rightarrow \infty}D_pY(u,n)D_q=P_{p}\,\eta(u)P_{q}=\zeta_{p,q}(u)
$$
where convergence is understood as described in Remark 3.1, which completes the proof.

\hfill $\blacksquare$
\begin{Corollary}
Under the assumptions of Theorem 4.2, it holds that 
$$
\lim_{n\rightarrow \infty}T_{p,q}(u,n)=
\eta_{p,q}(u)
$$ 
for any $p,q,u$
\end{Corollary}
{\it Proof.}
This is an easy consequence of Theorem 4.2.
\hfill $\blacksquare$

\begin{Example}
{\rm 
Computations of limit mixed (*-) moments of blocks reduce to the computation of moments
involving matricial $R$-circular systems of operators. For instance,
\begin{eqnarray*}
\lim_{n\rightarrow \infty}\tau_{q}(n)(S_{p,q}(n)^*S_{q,p}(n)^*S_{q,p}(n)S_{p,q}(n))&=&
\Psi_q(\zeta_{p,q}^*\zeta_{q,p}^*\zeta_{q,p}\zeta_{p,q})\\
&=&
\varphi(\ell_1^*\ell_1)\varphi(\ell_3^*\ell_3)\\
&=&d_pd_qv_{p,q}v_{q,p}
\end{eqnarray*}
where $p\neq q$, since 
$$
\zeta_{p,q}=(\ell_1+\ell_2^*)\otimes e(p,q), \;\;\;\zeta_{q,p}=(\ell_3+\ell_4^*)\otimes e(q,p)
$$
where $\{\ell_1, \ell_2, \ell_3, \ell_4\}$ is a *-free system of free creation operators with covariances
$\varphi(\ell_1^*\ell_1)=d_{p}v_{p,q}$ and $\varphi(\ell_3^*\ell_3)=d_{q}v_{q,p}$. In fact, it is well known that 
any pair of free circular operators can be written in the form $\ell_1 +\ell_2^*$ and $\ell_3+\ell_4^*$, respectively.}
\end{Example}

\section{Products of independent Gaussian random matrices}

The first concrete application of our method concerns products of independent rectangular Gaussian random matrices 
[10]. For any given $p\in {\mathbb N}$ and any $n\in {\mathbb N}$, consider the product of independent 
rectangular Gaussian random matrices 
$$
B(n)=X_1(n)X_2(n)\ldots X_{p}(n),
$$
where $n\in {\mathbb N}$ and all entries of these matrices are assumed to be independent $(0,1/n)$-Gaussian 
variables. If $X_{j}(n)$ is an $n_{j-1} \times n_{j}$ matrix for any $1\leq j\leq p$, we assume that 
$$
\lim_{n\rightarrow \infty} \frac{n_{j}}{n} =d_j
$$ 
for any $j\in \{0, \ldots , p\}$ (it is convenient to start with $n_{0}$ rather than with $n_1$).
Let $\tau_{0}(n)$ be the trace over the set of first $n_{0}$ basis vectors composed with classical expectation. 
\begin{Theorem}
Under the above assumptions, it holds that
$$
\lim_{n \rightarrow \infty}\tau_{0}(n)\left(\left(B(n)B^{*}(n)\right)^{k}\right)=P_{k}(d_0,d_1, \ldots , d_{p}),
$$
where
$$
P_{k}(d_0,d_{1},\ldots , d_{p})=\sum_{j_0 + \ldots + j_{p} = pk+1}\frac{1}{k}{k \choose j_0}{k \choose j_1}\ldots {k \choose j_{p}}
\;d_{0}^{j_{0}-1}d_{1}^{j_1}\ldots d_{p}^{j_{p}},
$$
for any natural $k$ and $j_0, j_1, \ldots , j_p$. These polynomials are called multivariate Fuss-Narayana polynomials. 
\end{Theorem}
{\it Sketch of the proof.}
The proof is based on embedding the matrices $X_1(n), \ldots , X_p(n)$ in symmetric blocks 
$T_{1,2}(n), \ldots , T_{p,p+1}(n)$, respectively, of a large square Gaussian random matrix $Y(n)$ 
of dimension $N\times N$, where $N=n_1+\ldots +n_{p+1}$, with $(0,1/n)$-Gaussian entries 
for any $n$ (we use only one matrix and thus we omit $u$ in our notations).
Computing the moments of $B(n)B^{*}(n)$ under $\tau_{0}(n)$ becomes now the normalized partial trace over the
first subset of basis vectors. We then use Corollary 4.1 to realize the limit moments in terms of 
$\eta_{1,2}, \ldots , \eta_{p,p+1}$ and their adjoints. These limit moments can be computed explicitly, which was done in [10], which completes the proof.
\hfill $\blacksquare$\\

Let us make some additional remarks on the above result:
\begin{enumerate}
\item
The special case of $p=1$ corresponds to Wishart matrices. If we set $d_{0}=1$ and $d_1=t$, we obtain 
$$
P_{k}(1,t)=\sum_{i+j=k+1} \frac{1}{k}{k \choose i}{k \choose j}t^{j},
$$
the moments of the Marchenko-Pastur distribution with
shape parameter $t>0$, namely
$$
\varrho_{t}={\rm max}\{1-t,0\}\delta_{0}+\frac{\sqrt{(x-a)(b-x)}}{2\pi x}1\!\!1_{[a,b]}(x)dx
$$
where $a=(1-\sqrt{t})^{2}$ and $b=(1+\sqrt{t})^{2}$. 
\item
The moments of the Marchenko-Pastur distributions are known to have the form
of {\it Narayana polynomials}
$$
N_{k}(t)=\sum_{j=1}^{k}\frac{1}{j}{k-1 \choose j-1}{k \choose j-1}t^{j}
$$ 
for any $k\in \mathbb{N}$. It is easy to show that $P_k(1,t)=N_{k}(t)$, but our formula is more suitable for multivariate generalizations.
\item
If $d_{0}=1$, the multivariate Fuss-Narayana polynomials become 
the moments of the $p$-fold free multiplicative convolution of the
Marchenko-Pastur distributions with arbitrary shape parameters. Namely, 
$$
P_{k}(1,t_{1},\ldots , t_{p})=m_{k}(\varrho_{t_1}\boxtimes \varrho_{t_2}\boxtimes \ldots \boxtimes \varrho_{t_{p}})
$$
where $m_k(\mu)$ stands for the $k$th moment of $\mu$ and $\boxtimes$ stands for the free multiplicative convolution. 
\item
A less general class of polynomials was studied in the context of {\it free Bessel laws} 
$$
\pi_{p,\,t}=\varrho_{1}^{\,\boxtimes (p-1)}\boxtimes \varrho_{t},
$$
where $p\in {\mathbb N}$, defined by Banica, Belinschi, Capitaine and Collins [2]. The moments of free Bessel laws are polynomials in one variable $t$, 
$$
Q_{k}(t)=\sum_{j=1}^{k}\frac{1}{j}{k-1 \choose j-1}{pk \choose j-1}t^{j}
$$
called {\it Fuss-Narayana polynomials}. 
Our polynomials are natural multivariate generalizations of these polynomials. Clearly, $Q_{k}(t)=P_{k}(1,\ldots , 1, t)$, where $1$ appears $p$ times.
\end{enumerate}

\section{New random matrix models}

The results on the asymptotics of random blocks can be applied to the construction of some new random matrix models, as we showed in [9].
In this paper, we choose to present a simple version of such a model for monotone independence [12]. 
In a similar way, one can construct 
random matrix models for boolean independence an s-free independence (freeness with subordination). 
In fact, more general versions, not restricted to Gaussian random matrices have been constructed in [8].

Monotone independence will appear in the study of the asymptotic joint distributions of  
two independent Hermitian Gaussian random matrices of the same block form. Our assumptions 
are the following:
\begin{enumerate}
\item[(A1)]
We have a family of independent Gaussian random matrices
$$
Y(u,n)=\left(
\begin{array}{rr}
A(u,n)& B(u,n)\\
C(u,n)& D(u,n)
\end{array}
\right)
$$
where $u\in \mathpzc{U}$, and blocks $(A(u,n))$ are evanescent, 
symmetric blocks built from $(B(u,n))$ and $(C(u,n))$ are unbalanced, and
blocks $(D(u,n))$ are balanced,
\item[(A2)]
the matrices are Hermitian, thus 
the off-diagonal blocks are Hermitian conjugate and the diagonal blocks are Hermitian
\item[(A3)]
the complex Gaussian variables $Y_{i,j}(u,n)$ have zero mean and have identical covariances within blocks, 
namely 
$$
\mathbb{E}(\overline{Y_{i,j}(u,n)}Y_{i,j}(u,n))=\frac{v_{p,q}(u)}{n}
$$ 
whenever the pair $(i,j)$ belongs to the block indexed by $(p,q)$,
\item[(A4)]
the decomposition of the identity matrix corresponding to the block decomposition is 
given by
$$
I(n)=D_1+D_2
$$
for any $n\in \mathbb{N}$.
\end{enumerate}

\begin{Theorem}
Under assumptions (A1)-(A4), if $\mathpzc{U}=\{1,2\}$ and $v_{p,q}(u)=1$ for any $p,q,u$, 
the pair $\{B(1,n)+C(1,n),Y(2,n)\}$ is asymptotically monotone independent with respect 
to the partial trace $\tau_1(n)$.
\end{Theorem}
{\it Proof.}
By Theorem 3.2, the proof reduces to showing that the pair 
$$
\{\omega_{2,1}(1), \omega_{2,1}(2)+\omega_{2,2}(2)\}
$$ 
is monotone independent with respect to $\Psi_1$ since, by assumption, 
the asymptotic dimensions are $d_1=0$ and $d_2=1$, which means that the remaining operators 
can be neglected. Denote $a=\omega_{2,1}(1)$ and $b=\omega_{2,1}(2)+\omega_{2,2}(2)$.
We need to show that 
$$
\Psi_1(w_1a_1b_1a_2w_2)=\Psi_1(b_1)\Psi_1(w_1a_1a_2w_2)
$$
for any $a_1,a_2\in {\mathbb C}[a,1_1]$, $b_1\in {\mathbb C}[b,1_2]$, where $1_1=1_{2,1}$ and $1_{2}=1$  and $w_1,w_2$ are arbitrary
elements of ${\mathbb C}\langle a,b,1_1,1_2\rangle$.
It suffices to consider the action of $a$ and $b$ onto their invariant subspace in ${\mathcal M}$
of the form
$$
{\mathcal M}'={\mathbb C}\Omega_1\oplus ({\mathcal F}(2)\otimes {\mathcal H}(1)) \oplus
({\mathcal F}(2)\otimes {\mathcal H}(2))
$$
where
${\mathcal F}(2)={\mathcal F}({\mathbb C}e_{2,2}(2))$ with the vacuum vector $\Omega$
and ${\mathcal H}(u)={\mathbb C}e_{2,1}(u)$ for $u\in \{1,2\}$, where we identify
$\Omega\otimes e_{2,1}(u)$ with $e_{2,1}(u)$. Now, the range of any polynomial in $a$ is 
${\mathbb C}\Omega_1\oplus {\mathcal H}(1)$ since
$$
a^{k}\Omega_1=
\left\{
\begin{array}{ll}
\Omega_1 & {\rm if}\;k\;{\rm is}\;{\rm even}\\
e_{2,1}(1) & {\rm if}\;k \;{\rm is}\;{\rm odd}
\end{array}
\right.
$$
and $1_1\Omega_1=\Omega_1$, $1_1e_{2,1}(1)=e_{2,1}(1)$. Therefore, it suffices to compute the action of
any polynomial in $b$ onto $\Omega_1$ and $e_{2,1}(1)$.  Now, the action of powers of $b$ onto $\Omega_1$ and onto $e_{2,1}(1)$
is the same as the action of the free Gaussian operator onto the vacuum vector in the free Fock space. Namely, we have
$$
b^{2k}\Omega_1=C_{k}\Omega_1\;{\rm mod}\;({\mathcal M}'\ominus {\mathbb C}\Omega_1)
$$
and
$$
b^{2k}e_{2,1}(1)=C_{k}e_{2,1}(1)\;{\rm mod}\;({\mathcal M}'\ominus ({\mathbb C}\Omega_1\oplus {\mathcal H}(1))
$$
where $C_k$ is the $k$th Catalan number and 
$$
b^{2k-1}\Omega_1= \;b^{2k-1}e_{2,1}(1)=0 \;{\rm mod}\; ({\mathcal M}'\ominus ({\mathbb C}\Omega_1\oplus {\mathcal H}(1))
$$ 
for any $k\in \mathbb{N}$. Thus $\Psi_1(b^{2k})=C_k$ and, moreover, since 
${\mathcal M}'\ominus ({\mathbb C}\Omega_1\oplus {\mathcal H}(1))\subset {\rm Ker}a$, 
the required condition for monotone independence holds if $b_1$ is a positive power of $b$. It is easy to see that it also
holds if $b_1=1_2$, which completes the proof.
\hfill $\blacksquare$\\

Interestingly enough, unbalanced symmetric blocks were also used in the construction of a 
simple random matrix model for free Meixner laws [1,4]. These are probability measures on the real line
associated with the sequences of Jacobi parameters of the form
$$
(\alpha_1, \alpha_2, \alpha_2, \ldots ) \;\;{\rm and}\;\;
(\beta_1, \beta_2, \beta_2, \ldots ) 
$$
and thus we can asy that they are associated with quadruples 
$(\alpha_1, \alpha_2, \beta_1, \beta_2)$. Let us formulate the theorem in the 
most interesting case when $\beta_1$ and $\beta_2$ are positive.
\begin{Theorem}
Under assumptions (A1)-(A4), let $\beta_1(u)=v_{2,1}(u)>0$ and $\beta_2(u)=v_{2,2}>0$, for any
$u\in \mathpzc{U}$. Then 
\begin{enumerate}
\item
the asymptotic distributions of the matrices
$$
M(u,n)=Y(u,n)+\alpha_1(u)D_1(n)+\alpha_2(u)D_2(n)
$$
under the partial trace $\tau_1(n)$ are the free Meixner distributions associated with the 
parameters $(\alpha_1(u), \alpha_2(u), \beta_1(u), \beta_2(u))$, respectively,
\item
the family $\{M(u,n):u\in \mathpzc{U}\}$ is asymptotically conditionally free with respect to the 
pair of partial traces $(\tau_1(n), \tau_2(n))$ as $n\rightarrow \infty$.
\end{enumerate}
\end{Theorem}
{\it Sketch of the proof.}
It follows from Theorem 3.2 that 
$$
\lim_{n\rightarrow \infty}Y(u,n)=\omega_{2,1}(u)+\omega_{2,2}(u)
$$
where the variances of $\omega_{2,1}(u)$ and $\omega_{2,2}(u)$ are $b_{2,1}(u)=d_2v_{2,1}(u)=\beta_{1}(u)$
and $b_{2,2}(u)=d_2v_{2,2}(u)=\beta_{2}(u)$ since $d_1=0$ and $d_2=1$. The fact that 
$d_1=0$ is the reason why $\omega_{1,1}(u)$ and $\omega_{1,2}(u)$ 
become trivial operators and that is why they do not show up in the limit realization.
Therefore,
$$
\lim_{n\rightarrow \infty}M(u,n)=\gamma(u):=\omega_{2,1}(u)+\omega_{2,2}(u)+\alpha_{1}(u)P_1+\alpha_{2}(u)P_2,
$$
where $P_1$ and $P_2$ are as in the proof of Theorem 3.2.
This implies, in particular, that  the asymptotic distribution of $M(u,n)$ under $\tau_{1}(n)$ agrees
with that of $\gamma(u)$ under $\Psi_1$. The proof that the moments of $\gamma(u)$ are the moments of the free Meixner law associated with the parameters $(\alpha_1(u), \alpha_2(u), \beta_1(u), \beta_2(u))$ 
is purely combinatorial and can be found in [9]. We also refer the reader to [9] for the proof of asymptotic conditional independence of the family of such matrices.
\hfill $\blacksquare$.

\end{document}